\documentclass[12pt]{amsart}
\usepackage{amssymb}
\usepackage{latexsym}
\usepackage{epsf}
\newtheorem{thm}{Theorem}[section]
\newtheorem{Thm}{Theorem}
\newtheorem{pro}[thm]{Proposition}
\newtheorem{lem}[thm]{Lemma}

\newtheorem{cor}[thm]{Corollary}
\newtheorem{Cor}[Thm]{Corollary}

\theoremstyle{definition}

\def\1{{\rm1\mathchoice{\kern-0.25em}{\kern-0.25em}
        {\kern-0.2em}{\kern-0.2em}I}}

\newcommand{\lmn}[1]{\vadjust{\setbox1=\vtop{\hsize 25mm
\parindent=0pt\baselineskip=9pt
\rightskip=4mm plus 4mm#1}
\hbox{\kern-26mm\smash{\raise .5ex\box1}}}}

\newcommand{\nc}{\newcommand}
\def\be#1\ee{\begin{equation}#1\end{equation}}
\nc{\bc}{\begin{center}} \nc{\ec}{\end{center}} \nc{\bb}{\mathbb}
\nc{\cal}{\mathcal} \nc{\frk}{\mathfrak} \nc{\N}{{\mathsf N}}
\nc{\K}{{\mathsf K}} \nc{\fk}{\mathbf{k}} \nc{\fn}{\mathbf{n}}
\nc{\fb}{\mathbf{b}}  \nc{\e}{\varepsilon} \nc{\ev}{{\rm{ev}}}

\theoremstyle{remark}

\def\Z{{\mathbb Z}}
\def\Q{{\mathbb Q}}
\def\N{{\mathbb N}}

\def\v8{\vskip 8pt}
\def\a{\alpha}

\def\la{\langle}
\def\ra{\rangle}
\def\l{\lambda}
\def\n{\nu}
\def\g{\gamma}
\def\m{\mu}

\def\fg{\mathfrak g}

\def\ZZ{\widehat{\Z[v]}_2}
\def\ZZZ{\widehat{\Z[v]}_s}
\def\Habiro{\widehat{\Z[q]}}

\begin{document}

\title[Laplace transform and universal invariants ]{Laplace transform and
universal $sl_2$ invariants}
\author{Anna Beliakova}
\address{Institut f\"ur Mathematik, Universit\"at Z\"urich,
 Winterthurerstrasse 190,
CH-8057 Z\"urich, Switzerland}
\email{anna@math.unizh.ch}
\author{Christian Blanchet}
\address{L.M.A.M., Universit\'e de Bretagne-Sud,
Centre de Recherche Tohannic, BP 573, F-56017 Vannes, France  }
\email{Christian.Blanchet@univ-ubs.fr}
\author{Thang Le}
\address{School of Mathematics, Georgia Institute of Technology, Atlanta,
GA 30332-0160, USA}
\email{letu@math.galtech.edu}
\date{September 2005}
\keywords{Universal quantum invariants, Habiro theory, Ohtsuki series,
cyclotomic completion ring}

\begin{abstract}
We develop a Laplace transform method
 for constructing   universal invariants
of $3$--manifolds.
As an application, we recover Habiro's theory
of integer homology $3$--spheres and extend it to some classes
of rational homology $3$--spheres with cyclic homology.
 If $|H_1|=2$,  we give explicit formulas
for universal invariants dominating
the  $sl_2$ and $SO(3)$ Witten--Reshetikhin--Turaev
 invariants, as well as their
 spin and cohomological refinements at {\it all} roots of unity.
New results on the Ohtsuki series and the integrality  of quantum invariants
 are the main applications of our construction.

\end{abstract}

\maketitle
\section*{Introduction}

For a simple Lie algebra $\fg$, the quantum invariants  of
$3$--manifolds (the Witten--Reshetikhin--Turaev invariants, see e.g.
\cite{Tu}) are defined only when the quantum parameter $q$ is a
certain root of unity. Habiro \cite{Ha} proposed a construction of a
universal $sl_2$ invariant of integer homology $3$--spheres
($\Z$HS), dominating all the quantum invariants. The results have
many important applications, among them are the integrality of
quantum invariants at {\em all} roots of unity, the recovery of
quantum invariants from the LMO invariant, and the possible
applications to the integral Topological Quantum Field Theory. All
the results were later extended to all simple Lie algebras by Habiro
and the third author \cite{HL}.

In this paper we extend Habiro's theory to some classes of rational
homology $3$--spheres and refined quantum invariants -- invariants of
spin structures and cohomological classes.

\subsection{Universal quantum invariants}
The universal invariant of a $\Z$HS is an element of the Habiro
ring
$$ \Habiro:=\lim_{\overleftarrow{\hspace{2mm}n\hspace{2mm}}}
\frac{ \Z[q]}{
(1-q)(1-q^2)...(1-q^n)}\, .$$
Every element $f\in \Habiro$ can be written as an infinite sum
$$ f(q)= \sum_{k\ge 0} f_k(q)\, (1-q)(1-q^2)...(1-q^n),$$
with $f_k(q)\in \Z[q]$. If $\xi$ is a root of unity, then $f(\xi)$ is
well--defined, since the summands become zero if $k$ is bigger than the
order of $\xi$. The Habiro ring has remarkable properties and is very
suitable for the study of quantum invariants. The result of Habiro
and Habiro--Le mentioned above is

\begin{Thm}{\rm (Habiro, Habiro--Le)}
For every simple Lie algebra $\fg$ and an integral homology $3$--sphere
$M$, there exists an invariant $I^\fg_M(q)\in \Habiro$, such that if
$\xi$ is a root of unity, then $I^\fg_M(\xi)$ is the quantum
invariant at $\xi$.
\end{Thm}

\subsubsection{Applications}

Let us mention the most important consequences of the Habiro's
construction. First of all, each product $(1-q)(1-q^2)\dots (1-q^n)$
is divisible by $(1-q)^n$, hence it is easy to expand every $f(q)\in
\Habiro$ into formal power series in $(q-1)$, denoted by $T(f)$ and
called the Taylor series of $f(q)$ at $q=1$. One important property
of $\Habiro$ is that $f\in \Habiro$ is uniquely determined by its
Taylor series. In other words, the map $T: \Habiro \to \Z[[q-1]]$
is injective. In particular, $\Habiro$ is an integral domain.
Another important property is that every $f\in \Habiro$ is
determined by the values of $f$ at any infinite set of roots of
unity of prime power  order. From the existence of $I^\fg_M$ one can
derive the following consequences  for $\Z$HS:
\begin{itemize}
\item
The quantum invariants   at all roots of unity are algebraic
integers.

\item The quantum invariants  at any infinite set of roots
 of unity of
prime power  order determine the whole set of quantum invariants.

\item
 Ohtsuki series (see \cite{Oh, Le2}) have integer coefficients and
 determines the whole set of quantum invariants.

 \item
 The Le--Murakami--Ohtsuki invariant (see \cite{LMO}) totally
 determines the quantum invariants.
\end{itemize}

The integrality of quantum invariants was established earlier only
at roots of unity of prime order (see \cite{MR,Le1}). The integrality
of the Ohtsuki series for $\fg=sl_2$ was proved by Rozansky, using
quite a different method.

\subsection{Results} In this paper we extend  Habiro's construction to
some classes of  rational homology $3$--spheres with cyclic
$H_1(M,\Z)$ and with
 spin/cohomological structures. Our results show that one can expect
  to generalize Habiro's theory to rational homology
spheres. For the case $H_1(M,\Z)=\Z/2\Z$, we get  the fullest results,
when all aspects of Habiro's theory are generalized.

We will use a new construction of the universal invariant based on
the {\em Laplace transform.}
This method originates from the paper of the third author \cite{Le}.
For $\Z$HS, this method reproduces  Habiro's results.

\def\cM{\mathcal M}
\def\RP{{\mathbold R}P}
\subsubsection{The case $H_1(M,\Z)=\Z/2\Z$} Let $\cM_b$ be the set of all
oriented closed 3-manifolds $M$ with
$H_1(M,\Z)=\Z/b\Z$.
If $M \in \cM_2$,  the quantum invariant $\tau_M$
depends on a square root $v$ of $q$.

When $q$ is an even root of unity, then the order of $v$
is divisible by 4. In this case, we put
$\tau'_M= \tau_M/\tau_{L(2,1)}$, i.e. we renormalize
the $sl_2$ quantum  invariant
  to be 1 for the lense space $L(2,1)$.

When $q$ is an odd root of unity, then $\tau_{L(2,1)}=0$ but the refined
$SO(3)$ version $\tau^{SO(3)}_{L(2,1)}\neq 0$.
We choose $v$  to be the root of $q$
which has the same order as $q$ does (i.e, also of odd order).
Here we put
$\tau'_M= \tau^{SO(3)}_M/\tau^{SO(3)}_{L(2,1)}$.

The role of $\Habiro$ will be replaced by
$$\ZZ := \lim_{\overleftarrow{\hspace{2mm}n\hspace{2mm}}}
\frac{\Z[v^{\pm 1}]}{(-v^2;-v)_{2n}},$$
where
$$ (-v^2;-v)_{2n} := \prod_{i=2}^{2n+1} (1+(-v)^i)= (1-v^3)(1-v^5)\dots (1-v^{2n+1}) \times
(1+q)(1+q^2) \dots (1+q^{n}).$$
Every $f(v) \in \ZZ$ can be written as, with $f_n(v) \in \Z[v^{\pm
1}]$,
$$f(v) = \sum_{n=0}^\infty f_n(v) \, (-v^2;-v)_{2n},$$
If
\begin{equation} \text{$v$ is a root of unity of order either odd or
divisible by 4}
\label{order}\end{equation} then $f(v)$ is well--defined. For every
root $q$ of unity, one can choose a square root $v$ of $q$ satisfying
(\ref{order}). The first main result is
\begin{Thm}\label{Two}
For every  closed oriented manifold $M\in \cM_2$, there exists an
invariant $I_M(v) \in \ZZ$, such that if $v$ is a root of unity
satisfying (\ref{order}), then $I_M(v)=\tau'_M(v)$.
\end{Thm}


Note that  $\ZZ$
embeds in $\Z[[v-1]]$, via Taylor series. As a consequence, we
will prove
\begin{Cor} \label{co2}
For $M\in \cM_2$ and the quantum invariants normalized so that
the projective space takes value 1, one has

(a) The quantum invariants
at all roots  of
unity are algebraic integers.

(b) The quantum  invariants at any infinite set $\{v\}$ of roots
of unity of odd
prime power  order determine the whole set of
quantum invariants.

(c) The
 Ohtsuki series, a formal power series in $q-1$, has  coefficients in $\Z[1/2]$, since it is equal to
 a formal power series in $v$, and
 determines the whole set of quantum invariants. If $v$ is a root of
 unity of order  $p^d$ with $p$ an odd prime, then the Ohtsuki
 series at $v$ converges $p$--adically to the quantum invariant at
 $v$.

 (d)
 The Le--Murakami--Ohtsuki invariant determines the quantum invariants at odd roots of 1.
\end{Cor}


\subsubsection{Spin structure and cohomological classes}
Suppose that the order  of $v$ is divisible by $4$, i.e. the order
of $q$ is even.  There are refined quantum invariants
$\tau_{M,\sigma}$, defined in \cite{KM}, where $\sigma$ is a spin
structure or a cohomological class in $H^1(M,\Z/2)$, depending on
whether the order of $v$ is equal to $0 \pmod 8$ or $4 \pmod 8$. We
will renormalize $\tau_{M,\sigma}$ by dividing by the
non--refined invariant of
the projective space, i.e $\tau'_{M,\sigma} : =\tau_{M,\sigma}/
\tau_{L(2,1)}$.
 Then we have
$\tau'_M=\sum_\sigma \tau'_{M,\sigma}$.
Let
$$\ZZZ := \lim_{\overleftarrow{\hspace{2mm}n\hspace{2mm}}}
\frac{\Z[v]}{(1+q)(1+q^2)\dots(1+q^n)}.$$ If $v$ is a root of
unity of order divisible by $4$, then $f(v)$ is well--defined for
$f\in \ZZZ$. For fixed $k$, if $n\geq 2k$ then
$(1+q)(1+q^2)\dots(1+q^n)$ is divisible by $(1+q)^k$, hence there
is a natural map, the Taylor series at $q=-1$, sending $f\in \ZZZ$
to $T_{-1}f \in \Z[I][1/2][[q+1]]$, where $I$ is the unit complex number.
 Habiro's theory \cite{Ha1}
shows that the map $T_{-1}$ is an embedding.
\begin{Thm}\label{main_spin}
For $M\in \cM_2$ and a spin structure (respectively, a cohomological
class) $\sigma$, there exists an invariant $I_{M,\sigma}(v) \in
\frac{1}{1-v}\ZZZ$, such that if $v$ is a root of unity of order divisible by 8
(respectively, equal to $4 \pmod 8)$, then $I_{M,\sigma}(v)$ is the
quantum invariant of $(M,\sigma)$ at $v$.
\end{Thm}

The integrality of $\tau_{M,\sigma}$ for $\Z/p\Z$--homology
spheres at roots of order  $2p$, where $p$ is an odd {\em prime}
and
 $\sigma$ is  a cohomological class,  was studied
by Murakami in \cite{HM, HM1}. Theorem \ref{main_spin} shows that
$(1-v) I_{M,\sigma}(v)$ is always an algebraic integer {\em for
all odd $p$} and that the quantum invariant $I_{M,\sigma}$ has an
expansion as formal power series in $(1+q)$. (The factor $(1-v)$
appears because we use the normalization for which the projective
space takes value 1).
 Theorems \ref{main_spin}, \ref{six} 
 give  partial answers to  Conjecture 5.3
and Remark 5.2  in \cite{HM}.
\v8

\noindent
{\bf Examples.} Suppose $M$ is obtained by surgery on the figure 8
knot with framing 2.
Then
$$ I_M(v) = \sum^\infty_{n=0} v^{-n(n+2)}(-v^2;-v)_{2n}$$

Suppose that the order of $v$ is $0 \pmod 8$, and this order
divided by $8$ is $\zeta\pmod 2$.
Let $\sigma_0$ be the characteristic spin structure on $M$,
and $\sigma_1$ the other one.
Then
$$ I_{M,\sigma_\e}(v)=\frac{1}{2(1-v)}
 \sum^\infty_{n=0} v^{-n(n+2)}\prod^n_{i=1}(1+q^i)
\left[\prod^n_{i=0} (1-v^{2i+1})- (-1)^{\zeta+\e}\prod^n_{i=0}(1+v^{2i+1})\right]
\, .$$

Assume that $v$ is $4p$--th root of unity with odd $p$ and
$v^{p^2}=\zeta I$, where $I$ is
the unit complex number and $\zeta=\pm 1$.
Let  $\sigma_\e\in H^1(M,\Z/2\Z)$, and $\sigma_1$ be trivial.
Then
$$ I_{M,\sigma_\e}(v)=\frac{1}{2(1-v)}
 \sum^\infty_{n=0} v^{-n(n+2)}\prod^n_{i=1}(1+q^i)
\left[\prod^n_{i=0} (1-v^{2i+1})+
(-1)^{\e}\zeta I\prod^n_{i=0}(1+v^{2i+1})\right] \, .$$

\subsubsection{The case  $H_1(M,\Z)=\Z/b\Z$}
Let $M\in \cM_b$. Assume that
the greatest common divisor of $b$ and the order $r$ of
 $q$ is a power of two.
More precisely,
 we suppose that $b=2^t c$ and $r=2^s d$
with odd $c,d$ and $gcd(c,d)=1$.
If $r$ is even and $t\neq s+1$, then
$\tau_{L(b,1)}\neq 0$ and we can renormalize
$\tau'_M=\tau_M/\tau_{L(b,1)}$ and $\tau'_{M,\sigma}=\tau_{M,\sigma}/\tau_{L(b,1)}$.
For odd $r$ ($s=0$), we put
$\tau'_M=\tau^{SO(3)}_M/\tau^{SO(3)}_{L(b,1)}$.
We show that the Laplace transform method works  and leads to formulas for
universal quantum invariants and  their refinements.
As a consequence, we have

\begin{Thm}\label{six}
Let $b=2^tc$ with odd $c$. Let $S=\{2^s d\in \N: gcd(c,d)=1,
d\;{\rm odd}, s\neq t-1 \}$. Let $M\in\cM_b$.
The  quantum invariants $\tau'_{M}$
 at   roots of unity of order
 $r\in S$  are algebraic integers. If $t>1$, then also
the refined quantum invariants  $\tau'_{M,\sigma}$ at even  roots of unity
of order $r\in S$  are algebraic integers.
\end{Thm}

\subsection{Plan of the paper}
The paper is organized as follows. After introducing
the Laplace transform method, we apply it to $\Z$HS
and get precise formulas for  Habiro's
universal invariants.
Then we apply this method to $\Q$HS with $|H_1|=2$.
Here again the exact formula for the  Laplace transform
 implies  various above mentioned results.
After that, refinements of quantum invariants are considered.
In Section 4, we derive
 explicit formulas for the spin and cohomological
refinements of universal invariants assuming  $|H_1|=2$.
In Section 5,
 we
construct  refined universal invariants in the case when
$H_1=\Z/b\Z$
and the greatest common divisor of $r$ and $b$ is a power of two.

\subsection*{Acknowledgment} The first author wishes to express her
gratitude to Dennis Stanton for the significant simplification of the proof
of Lemmas \ref{l1}, \ref{l2}.

\section{Laplace transform}

In this section we introduce the Laplace transform method.
\subsection{Cyclotomic expansion of the colored Jones polynomial}
Let $K$ be a knot with framing zero.
 We denote by $J'_K(\l)$ the  Jones
polynomial of $K$ colored by the $\l$--dimensional irreducible
representation of $\mathfrak sl_2$, and  normalized at one for the
unknot. Note that $J'_K(\l)\in \Z[q^{\pm 1}]$.

In \cite{Ha}, Habiro announced that there exist $C_{K,k}\in \Z[q^{\pm 1}]$
such that \be\label{hab}
 J'_K(\l)=\sum^{\infty}_{k=0} C_{K,k}\; (q^{1+\l})_k (q^{1-\l})_k\, .\ee
Here we use the standard notation
$(a)_n=(1-a)(1-aq)(1-aq^2)...(1-aq^{n-1})$. The sum in (\ref{hab})
is finite, because the summands with $k\geq \l$ are zero. This
expansion is called the cyclotomic expansion of the colored Jones
polynomial. The nontrivial part here is that $C_{K,k}$'s are Laurent
polynomials in $q$ with integer coefficients.
\v8

\noindent {\bf Examples.} For the right--, left--handed trefoil and
the figure 8 knot,  we have
$$J'_{3_1}(\l)=\sum^{\infty}_{k=0} q^{-{k(k+2)}}(q^{1+\l})_k (q^{1-\l})_k \,
$$
$$J'_{\bar 3_1}(\l)=\sum^{\infty}_{k=0} q^{k}(q^{1+\l})_k (q^{1-\l})_k \,
$$
$$J'_{4_1}(\l)=\sum^{\infty}_{k=0} (-1)^k q^{-\frac{k(k+1)}{2}}(q^{1+\l})_k (q^{1-\l})_k \,
.$$
\noindent
{\bf Note.} The coefficients $C_{K,k}$ are computed for all twist knots
in \cite{Ma}.

\subsection{Quantum invariants for knot surgeries}
Let $M=S^3(K_b)$ be a $\Q$HS obtained by surgery on $K$ with
nonzero framing $b$. Assume that $q$ is
 a primitive $r$--th root of unity and $r$ is even.
The  quantum $sl_2$ invariant of $M$ is defined as follows, see
\cite{Tu}.
 \be\label{inv} \tau_M(q)= \frac{\sum\limits^{r-1}_{
{\l=0}}\; q^{\frac{b(\l^2-1)}{4}} \;\,
(1-q^\l)(1-q^{-\l})J'_K(\l)}{{\sum\limits^{r-1}_{\l=0\; }
q^{\frac{sn(b)(\l^2-1)}{4}}\;\, (1-q^\l)(1-q^{-\l})}\,} , \ee where
$sn(b)$ is the sign of $b$. To be precise, one needs to fix a 4--th
root of $q$. Note that when computing the Jones polynomial
of a knot (or a link) in this paper, we always assume that its framing
is zero. However in the formula for the quantum invariant,
 framing is taken into account by means of the factor $q^{b(\l^2-1)/4}$.

Substituting Habiro's formula (\ref{hab}) into (\ref{inv}) we get
\be\label{inv1} \tau_M(q)= \frac{\sum\limits^{r-1}_{\l=0}\;
q^{\frac{b(\l^2-1)}{4}} \;\,\sum\limits^\infty_{n=0}\,
  C_{K,n} F_n(q^\l,q)}{{\sum\limits^{r-1}_{\l=0} }\;
q^{\frac{sn(b)(\l^2-1)}{4}}\;\, F_0(q^\l,q)}\, , \ee where
$F_n(q^\l,q)=(q^\l)_{n+1}(q^{-\l})_{n+1}$.

Suppose  $r$ is odd. Then, taking the sum over odd $\l$ in the
numerator and the denominator of (\ref{inv}) we get the  $SO(3)$
invariant of $M$. In this case, there is no need to fix  4--th root of
$q$.


\subsection{Laplace transform method}

The main idea behind the Laplace transform method is to interchange
the sums over $\l$ and $n$ in (\ref{inv1}) and regard
$\sum^{r-1}_{\l } q^{b(\l^2-1)/4}$ as an operator (called Laplace
transform) acting on  $F_n(q^\l,q)$.

More precisely, after interchanging the sums in the numerator of
(\ref{inv1}) we get
$$\sum^{r-1}_{n=0}C_{K,n}(q)\sum^{r-1}_{\l }
q^{\frac{b(\l^2-1)}{4}} F_n(q^\l,q)\, .$$ Now  observe, that
$F_n(q^\l,q)=(q^\l)_{n+1}(q^{-\l})_{n+1}$ is a polynomial in two
variables $q^\l$ and $q$. The Laplace transform does not affect $q$,
and  we only need  to compute the action of Laplace on  $q^{a\l}$.

Suppose the greatest common divisor of $b$ and $r$ is 1 or 2,
and $r$ is even. A
simple square completion argument shows that
$$\sum^{r-1}_{\l=0} q^{\frac{b(\l^2-1)}{4}}\, q^{a\l}=
q^{-\frac{a^2\, b^*}{gcd(b,r)}}\,
\gamma_{b,r} \, $$ where $b^*$ is an integer such
that $b^* b =gcd(b,r) \pmod r$, and
$$\gamma_{b,r}:=\sum^{r-1}_{\l=0} q^{\frac{b(\l^2-1)}{4}}\, . $$
Summarizing the previous discussion, we get
$$\sum^{r-1}_{\l =0}
q^{\frac{b(\l^2-1)}{4}}\, F_n(q^\l,q)= \ev_r(L_b(F_n(q^\l,q)))\;
\gamma_{b,r}\, .$$ Here $L_b(F)$ is the Laplace transform of $F$,
which is defined as follows. Suppose  $F$ is a formal power series
in $q^{\pm 1}$ and $q^{\pm \l}$. Then $L_b(F)$ is obtained from $F$
by  replacing every $q^{a\l}$  by $q^{-a^2/b}$. The evaluation map
$\ev_r$ converts $q^{1/b}$ to $(q^{1/gcd(b,r)})^{b*}$. Note that
while $\ev_r$ might depend on $r$, the Laplace transform $L_b$ does
not. And also if $b=1$ or $b=2$, then $\ev_r$ does not depend on
$r$: In these cases, $\ev_r(q^{1/b})=q^{1/b}$.

If $r$ is odd and $gcd(b,r)=1,2$, we can define the Laplace transform
by the same formula (i.e. $q^{a\l} \mapsto q^{-a^2/b}$).
 In this case, we have
$$\sum^{r-2}_{\l =1 \; {\rm odd}}
q^{\frac{b(\l^2-1)}{4}}\, F_n(q^\l,q)= \ev_r(L_b(F_n(q^\l,q)))\;
\gamma^1_{b,r}\, ,$$  where
$$\gamma^1_{b,r}:=\sum^{r-2}_{\l=1 \; {\rm odd}} q^{\frac{b(\l^2-1)}{4}}\, . $$


As a result,  we have  closed formulas for  quantum
invariants in terms of the Laplace transform.

\begin{thm}\label{le}
Let $M=S^3(K_b)$ and $gcd(b,r)$ divide $2$.  Then
$$\tau_M(q)=\frac{1}{2(1-q^{-sn(b)})} \frac{\gamma_{b,r}}{\gamma_{sn(b),r}}\;
\sum^\infty_{n=0} C_{K,n} \ev_r(L_b(F_n))\, ,$$
$$\tau^{SO(3)}_M(q)
=\frac{1}{2(1-q^{-sn(b)})} \frac{\gamma^1_{b,r}}{\gamma^1_{sn(b),r}}\;
\sum^\infty_{n=0} C_{K,n} \ev_r(L_b(F_n))\, .$$

\end{thm}

\section{Habiro theory}
In this section we show how Theorem \ref{le}
can be used to compute Habiro's
 universal  invariants of $\Z$HS.

\subsection{Knot surgeries}
Any knot surgery with framing $b=\pm 1$ yields a $\Z$HS. Combining
Theorem \ref{le} with Lemma \ref{l1} below we get the following
theorem.

\begin{thm}\label{le1} (Habiro)
For $M_\pm=S^3(K_{\pm 1})$, we have
$$ \tau_{M_+}(q)=
\sum^\infty_{n=0} (-1)^{n}  q^{-\frac{n(n+3)}{2}} C_{K,n}(q)
\frac{(q^{n+1})_{n+1}}{1-q}\, ,$$
$$ \tau_{M_-}(q)=
\sum^\infty_{n=0} C_{K,n}(q) \frac{ (q^{n+1})_{n+1}}{1-q}\, .$$
\end{thm}

\noindent
{\bf Remark.}  The formulas in Theorem \ref{le1} do not depend on the
order of the root of unity $q$, and, in fact,
define
 elements of the Habiro's ring
which dominate  quantum invariants at all roots of unity and,
therefore, have to coincide with the Habiro's universal $sl_2$
invariants of $M_{\pm}$.

\v8

\noindent {\bf Examples.} Denote by $3_1$ and $4_1$ the Poincare
sphere and
 the $3$--manifold obtained by framing 1
surgery on figure 8 knot. By Theorem \ref{le1}, we have
$$\tau_{3_1}(q)= \frac{q}{1-q}\sum^\infty_{k=0}(-1)^k
q^{-\frac{(k+2)(3k+1)}{2}} (q^{k+1})_{k+1}$$
$$\tau_{4_1}(q)= \frac{q}{1-q}\sum^\infty_{k=0}(-1)^k
q^{-(k+1)^2} (q^{k+1})_{k+1}$$

\begin{lem} \label{l1}
\be\label{lap1}L_{-1}((q^\l)_{k+1}(q^{-\l})_{k+1})=2(q^{k+1})_{k+1}\, . \ee
\be\label{lap2}
L_1((q^\l)_{k+1}(q^{-\l})_{k+1})=2(-1)^{k+1} q^{-\frac{(k+2)(k+1)}{2}}
(q^{k+1})_{k+1}\, .\ee
\end{lem}

\begin{proof}
First, note that (\ref{lap2}) follows from (\ref{lap1})
and
$$L_{-b}(F_n(q^\l,q))= q^{k(k+1)}\, L_b(F_n(q^\l,q^{-1}))\, .$$
Let us prove (\ref{lap1}). For this, we split
$$F_k(q^\l,q)=S_k(q^\l,q)+T_k(q^\l,q)$$
with $S_k(q^\l,q)=(q^\l)_{k+1}(q^{-\l+1})_k$ and
$T_k(q^\l,q)=-q^{-\l} (q^\l)_{k+1}(q^{-\l+1})_k$.
Then $S_k(q^{-\l},q)=T_k(q^\l,q)$ implies
$L_{b}(S_k)=L_{b}(T_k)$
for any $b$. Therefore, we have to look at one of them only.

Further, by the $q$--binomial theorem (eq. (II.4) in \cite{GR}) we get
$$S_k(q^\l,q)=(-1)^k q^{-k\l}q^{k(k+1)/2}(q^{\l-k})_{2k+1}=$$
$$
(-1)^kq^{\frac{k(k+1)}{2}}\sum^{2k+1}_{j=0} (-1)^j\left[\begin{array}{c}
2k+1\\ j\end{array}\right]_q q^{\frac{j(j-1)}{2}}q^{-kj}q^{(j-k)\l}\, $$
where
$$\left[\begin{array}{c}n\\k\end{array}\right]_q
=\frac{(q)_n}{(q)_k(q)_{n-k}}\, .
$$
Taking the Laplace transform we have
$$L_{-1}(S_k(q^\l,q))=(-1)^k q^{\frac{3k^2+k}{2}}\sum^{2k+1}_{j=o}\frac{
(q^{-2k-1})_j}{(q)_j}q^{j^2+j-jk}\, .$$
The result follows now by applying the
 Sears--Carlitz transformation (eq. (III.14)
in \cite{GR}) for terminating $_3\phi_2$ series with specializations
$a=q^{-2k-1}$, $b,c\to \infty$, $z\to q^{k+2}$.
\end{proof}

\subsection{Link surgeries}\label{2.2}

Analogous to the case of knots,
Habiro gave an expression for
the colored Jones function of links.
To introduce his formula we need  some notation.

Let $L$ be an algebraically split framed link of $l$ components in
$S^3$ with all framings zero.
Let $\fn=\{n_1, n_2,...,n_{l}\}$ be a coloring of $L$ by
$\fn$--dimensional  irreducible
representations of ${\mathfrak sl}_2$.
 We denote by $J_L(\fn)$ the $\fn$--colored Jones polynomial
of $L$. We put
$$J'_L(\fn)=\frac{J_L(\fn)}{[\fn]}\, ,$$
where $[\fn]=\prod_i [n_i]$ with
 $[i]=(v^i-v^{-i})/(v-v^{-1})$, and $v^2=q$.
 Theorem 3.3 in \cite{Ha} implies then the following.
More details are given in Appendix.

\begin{pro}(Habiro) \label{link}
There exist $C_{L,\fk}(v)\in \Z[v^{\pm 1}]$ such that
$$J'_L(\fn)=\sum_{k=0}^\infty   \left(\sum_{\max k_i=k} C_{L,\fk}(v) \, (1-q)^l  \prod^{l}_{i=1}
\frac{(q^{1+n_i})_{k_i}(q^{1-n_i})_{k_i}}{(q^{k_i+1})_{k_i+1}}\right)
\frac{(q^{k+1})_{k+1}}{ (1-q)}
$$
\end{pro}
\v8

\noindent
{\bf Example.}
Let $L$ be the $0$--framed Whitehead link.
$$ J'_L(\l,\m)=\sum^{\infty}_{k=0} (-1)^{k} v^{-k(k+1)}(1-q)
\frac{(q^{1+\mu})_{k}(q^{1-\mu})_{k} }{(q^{k+1})_{k+1}}\;
(q^{1+\l})_{k}(q^{1-\l})_{k}$$
\v8

Let $M=S^3(L)$ be obtained by surgery on the framed link $L$ of $l$
 components in $S^3$.
 We denote by $b_i$ the framing of the $i$--th component
of $L$. Let $\sigma_+$ (respectively, $\sigma_-$) be the number of positive
(respectively, negative) eigenvalues of the linking matrix for $L$.
We put
$$Q_L(\fn):=  J_L(\fn)\times [\fn].$$
By definition, the quantum invariant of $M$ is

\begin{equation} \label{def} \tau_M(q) =
\frac{\sum_{\fn}^{\bf r-1}
 \prod^l_{i=1}  q^{b_i(n^2_i-1)/4}
 Q_L(\fn)}{(\sum_{n=0}^{r-1} q^{(n^2-1)/4}
 [n]^2)^{\sigma_+}\, (\sum_{n=0}^{r-1} q^{-(n^2-1)/4} [n]^2)^{\sigma_-}}
\end{equation}
\v8

Suppose  $M$ be a $\Z$HS.  Without loss of generality, we  can assume that
 $L$ is an algebraically split link with framings $\pm 1$.
Suppose that the  first $\sigma_+$ components have framing $+1$,
and the others $-1$. Substituting cyclotomic expansion
of the colored Jones polynomial (given in Proposition \ref{link})
into (\ref{def})
and applying the Laplace transform method to each component of $L$,
we derive the following formula for the
  universal $sl_2$ invariant of $M$.
\begin{thm} (Habiro) For $M$ as above, we have
\be\label{inv11} \tau_M(q)= \sum_{k=0}^\infty   \left(\sum_{\max
k_i=k} C_{L,\fk}(v)
 \prod\limits^{\sigma_+}_{i=1}
(-1)^{k_i} q^{-\frac{k_i(k_i+3)}{2}}\right) \frac{(q^{k+1})_{k+1}}{
(1-q)} \, . \ee

\end{thm}

Again, the right hand side belongs to $\Habiro$ and defines the
universal invariant of Habiro. Note that the $SO(3)$ invariant
of $M$ is also given by (\ref{inv11}).

\v8

\section{Rational homology $3$--spheres with $|H_1(M)|=2$}\label{so3}

In this section we define universal invariants of $\Q$HS with $|H_1|=2$.

\subsection{Normalization}
Suppose that the order of $v$ is divisible by 4.
The projective space, or the lense space $L(2,1)$ should be
considered as the unit in this class. It's easy to show that the
quantum invariant of $L(2,1)$ is given by
$$ \tau_{L(2,1)}(v) = \frac{\gamma_{2,r}}{(1+v^{-1})\,
\gamma_{1,r}} = \frac{\gamma_{-2,r}}{(1+v)\, \gamma_{-1,r}}.$$
For $M\in \cM_2$, we will use a normalization such that the
projective space $L(2,1)$ takes value 1:
$$ \tau'_M := \tau_M/\tau_{L(2,1)}.$$
If $v$ is an odd root of unity, we put
$$ \tau'_M := \tau^{SO(3)}_M/\tau^{SO(3)}_{L(2,1)},$$
where $$\tau^{SO(3)}_{L(2,1)}=\frac{\gamma^1_{2,r}}{
(1+v^{-1})\gamma^1_{1,r}}=\frac{\gamma^1_{-2,r}}{
(1+v)\gamma^1_{-1,r}}\, .$$

\subsection{Universal invariants}
Let $M_{\pm}=S^3(L)$, where $L$ is an $(l+1)$--component link numbered
by $0,1,\dots,l$. Assume that the $0$--th component has framing $\pm
2$, the next $s$ components have framing 1, and the remaining ones
have framing $-1$.

\begin{pro}\label{h2} For $M_\pm$ as above, we have
$$\tau'_{M_+}(v)=\sum_{k=0}^\infty   \left(\sum_{\max
k_i=k} C_{L,\fk}(v)\, (-v)^{-k_0}\prod_{i=k_0+1}^{k} (1+v^{2i+1})
 \prod\limits^{s}_{i=1}
(-1)^{k_i} q^{-\frac{k_i(k_i+3)}{2}}\right) (-v^2;-v)_{2k}
 \, .
$$

$$\tau'_{M_-}(v)=\sum_{k=0}^\infty   \left(\sum_{\max
k_i=k} C_{L,\fk}(v)\prod_{i=k_0+1}^{k} (1+v^{2i+1})
 \prod\limits^{s}_{i=1}
(-1)^{k_i} q^{-\frac{k_i(k_i+3)}{2}}\right) (-v^2;-v)_{2k}
 \, .
$$

\end{pro}

Note that $\tau'_{M_\pm}\in \ZZ$.
Theorem \ref{Two}
follows from
Proposition \ref{h2} and Lemma \ref{top} below, which states
that every $M\in \cM_2$ can be obtained from $S^3$ by surgery along a
link as described.

\v8
\noindent
{\bf Example.} Let $L$ be the Whitehead link with framings $2$ and $-1$.
Let $M=S^3(L)$.
$$\tau'_M(v)=\sum^\infty_k v^{-k(k+2)} (-v^2;-v)_{2k} $$
\v8

\begin{proof}
The proof is again an application of the Laplace transform method, and
Lemma \ref{l2} below. In addition, we use the
following identity
$$\frac{(-v^2;-v)_{2k_0}}{(q^{k_0+1})_{k_0+1}}\; (q^{k+1})_{k+1}=
(-v^2;-v)_{2k}\,  \prod^{k}_{i=k_0+1} (1+v^{2i+1}) \; ,$$ whose
proof is left to the reader.
Clearly, the formulas in Lemma \ref{l2} remain true after
replacing $\g_{b,r}$ with $\g^1_{b,r}$ and $\tau_{L(2,1)}(v)$
with $\tau^{SO(3)}_{L(2,1)}(v)$.

\end{proof}

\begin{lem}\label{l2}

$$\frac{L_{2}[(q^\l)_{k+1}(q^{-\l})_{k+1}]\, \gamma_{2,r}}{2(1-q^{-1})
\gamma_{1,r}}=(-v)^{-k}\, (-v^2;-v)_{2k}\, \tau_{L(2,1)}(v)$$
\be\label{b2} \frac{L_{-2}[(q^\l)_{k+1}(q^{-\l})_{k+1}]\,
\gamma_{-2,r}}{2(1-q)\, \gamma_{-1,r} }=(-v^2;-v)_{2k}\,
\tau_{L(2,1)}(v)\ee
\end{lem}

\begin{proof}

We proceed by proving (\ref{b2}).
By the $q$--binomial theorem we get
$$L_{-2}(F_k(q^\l,q))=
2 (-1)^kq^{k^2+k/2}\sum^{2k+1}_{j=0} \frac{(q^{-2k-1})_j}{(q)_j}
q^{j+j^2/2}\, .$$
The Sears--Carlitz transformation (eq. (III.14) in \cite{GR}) with
$a=q^{-2k-1}$, $c=-q^{-k}$, $z=q^{k+3/2}$ and $b\to \infty$
reduce this sum to $_2\phi_1(-q^{-k-1/2}, q^{-k};q^{-k+1/2},q)$
which can be computed by the $q$--Vandermode formula (eq. (II.6) in \cite{GR}).
As a result, we get
$$L_{-2}(F_k)=(1-v)(-v^2;-v)_{2k}\, ,\;\;\;\;\;
L_2(F_k)=(-1)^{k+1}v^{-k-1}(1-v)(-v^2;-v)_{2k}\, .$$

\end{proof}

\begin{lem}\label{top}
Any $M\in\cM_2$ can be obtained from $S^3$ by surgery on an
algebraically split link with framing $\pm 2$ on one component
and framings $\pm 1$ on the others.
\end{lem}
\begin{proof}
Choose a loop $K$ representing the nontrivial homology class
of $M$. Then $M\setminus K$ has homology of a solid torus.
By doing an integral surgery on $K$, we get a $\Z$HS $M'$.
In $M'$, $K$ spans a surface. We shrink the surface to its core,
i.e. a 1--dimensional complex. Now $M'$ can be obtained from $S^3$
by surgery on an algebraically split link $L$. Furthermore, $L$
  can be isotoped to miss the core of the spanning surface.
Hence, $L\cup K$ is an algebraically split surgery link for $M$
satisfying the required conditions.
\end{proof}
\v8

\noindent
{\bf Proof of Corollary \ref{co2}.}
 By  Theorem 5.4 in \cite{Ha1}, there exists an injective
homomorphism $\ZZ \to \Z[[1-v]]$ generating Ohtsuki series. More
details  will be given in \cite{BL}.

\section{Refinements}
In this section we show that the Laplace transform method can
effectively be used also to define refined universal invariants.

Suppose $\sigma$ is a spin structure (respectively, a cohomological
class in $H^1(M,\Z/2)$) and the order of $v$ is divisible by 8
(respectively, is equal to $4\pmod 8$), then there is defined the
refined invariant $\tau_{M,\sigma}(v)$.
We will
use the normalization
$$ \tau'_{M,\sigma}(v) =\frac{ \tau_{M,\sigma}(v)}{\tau_{L(2,1)}(v)}\, ,
\;\;{\rm  }\;\;\;\;
\tau'_M=\sum_\sigma \tau'_{M,\sigma}\, .
$$

\subsection{Spin refinements for  $M\in \cM_2$}

Without loss of generality, we
 will  assume that $M$ is obtained by surgery along the link $L$
of $(l+1)$ components,
as described in the previous section. The framing of the
$0$--th component is $\eta 2$, where $\eta=\pm 1$. Then $M$ has 2
spin structure $\sigma_0$ and $\sigma_1$, corresponding to the two
characteristic sublinks: one is the whole $L$ and the other is $L$
with the $0$--th component removed.

In this subsection
we suppose that $q$ is an $r$--th root of unity
 of order divisible by 4 and $v^2=q$.
By definition,
\be \tau_{M,\sigma_\e}(v) = \frac{\sum_{n_0\equiv \e
({\rm mod}\, 2)}^{r-1}
\sum_{n_1,n_2,\dots ,n_l\, \text{even}}^{r-1}
 q^{\eta(n_0^2-1)/2} \prod^s_{i=1}  q^{(n^2_i-1)/4}
\prod^l_{i=s+1}  q^{-(n^2_i-1)/4} Q_L(\fn)}{(\sum_{n=0}^{r-1} q^{(n^2-1)/4}
 [n]^2)^{s+\eta}\, (\sum_{n=0}^{r-1} q^{-(n^2-1)/4} [n]^2)^{l+1-s-\eta}}
 \label{def_spin}
 \ee
The next lemma  is well--known (compare \cite{KM}, \cite{B}).
 \begin{lem}\label{wk}
$$\sum_{n_1,n_2,\dots ,n_l
\text{even}}^{r-1}  \prod^s_{i=1}  q^{(n^2_i-1)/4}
\prod^l_{i=s+1}  q^{-(n^2_i-1)/4}\;
Q_L(\fn) =$$
$$ \sum_{n_1,n_2,\dots ,n_l
=0}^{r-1}  \prod^s_{i=1}  q^{(n^2_i-1)/4}
\prod^l_{i=s+1}  q^{-(n^2_i-1)/4} \; Q_L(\fn).$$
 \end{lem}

\begin{Thm}\label{four} Suppose the order $2r$ of $v$ is divisible by 8. Let $r/4\equiv \zeta \pmod 2$.
Then
\begin{equation} \tau'_{M,\sigma_\e}(v) =
\frac{1}{2}\left[ \tau'_M(v) -\eta (-1)^{\zeta +\e}
\frac{(1+v)}{(1-v)}\tilde\tau_M(-v) \right]\, ,\label{10}
\end{equation}
where $\tilde\tau_M (-v)$ is obtained from $\tau'_M(-v)$
given in Proposition \ref{h2}
by  replacing $C_{L,\fk}(-v)$ with $C_{L,\fk}(v)$.
\end{Thm}

The proof will be given in the next subsection.
 It's easy to
see that the right hand side of (\ref{10}) belongs to $\frac{1}{1-v}\ZZZ$,
and define an invariant of $3$--manifold $M\in \cM_2$ with a fixed
spin structure. This proves the part  of Theorem \ref{main_spin}
concerning the spin structure.

\subsection{Proof of Theorem \ref{four}}
Let us first introduce the odd and even Laplace transforms
as follows.
For $\e=0,1$, we put
$$\gamma^\e_{b,r}=\sum^{r-1}_{\l=\e \mod 2} q^{b(\l^2-1)/4}\,.  $$
We set

\begin{equation}
L^\e_b (P(q^\l,v)) := \frac{1}{ \gamma_{b,r}}  \, \sum^{r-1}_{
\l=\e \mod 2} q^{b(\l^2-1)/4}\, P(q^\l,v)\,
\label{new1}\end{equation} where $P(q^\l,v)$ is a Laurent
polynomial in $q^\lambda$ and $v$.

Let us prove Theorem \ref{four}
assuming Lemma \ref{as} below.
To compute the invariant, we  need to insert the cyclotomic
  expansion
of the colored Jones polynomial
  (given in Proposition \ref{link}) into (\ref{def_spin})
 and use the Laplace transform method.
By Lemma \ref{wk}, we need to apply $L_{\pm 1}$ to all components
except of the $0$--th one, and  $L^\e_{\pm 2}$ to the $0$--th component.
From $L^0_{\pm 2}+L^1_{\pm 2}=L_{\pm
  2}$
 and c) of the lemma we get
$$L^\e_{\pm 2}=\frac{1}{2}\left(\; L_{\pm 2}+
  (-1)^{\e+1}(c_1-c_0)L_{\pm 2}|_{v\to -v}
\;\right)\, . $$
The constants $c_\e$ are given in the proof of Lemma \ref{as}.
The result follows now from the next two formulas:
$$\frac{L_{2}(F_k)|_{v\to -v}\, }{2(1-q^{-1})}\;\frac{\gamma_{2,r}}
{\gamma_{1,r}}=-\frac{v^{-k}(1+v)}{1-v}\, (-v^2;v)_{2k}\, \tau_{L(2,1)}(v)$$
$$ \frac{L_{-2}(F_k)|_{v\to -v}\,
}{2(1-q)}\,\frac{\gamma_{-2,r}}{ \gamma_{-1,r} }=\frac{1+v}{1-v}(-v^2;v)_{2k}\,
\tau_{L(2,1)}(v)$$
\v8
$\hfill\Box$

\begin{lem}\label{as}
There are exist  constants $c_\e$, independent on $r$, such that
$c_0+c_1=1$, and

a) $\g^\e_{\pm2,r}=c_\e \g_{\pm2,r}$;

b) $L^\e_{\pm2} (q^{a\l})=c_{\e+a} L_{\pm2}(q^{a\l})$, where
$\e+a$ is taken modulo $2$;

c)
$(L^1_{\pm2}-L^0_{\pm2})(q^{a\l})=(c_1-c_0)L_{\pm2}(q^{a\l})|_{v\to
-v} \, .$

\end{lem}

\begin{proof}
a) By shifting $\l \to\l+r/2$, we see that
$\g^0=0$ if $r=4p$  ($p$ odd) or $\zeta=1$,
 and $\g^1=0$ if $r$ is divisible
by 8 or $\zeta=0$.
This implies $c_0=0$ in the first case, and $c_1=0$ in the second one.

b) If $r=4p$ ($p$ odd), we have
(compare with the case $s=t+1$ in the next section)
$$L^0_{\pm2}(q^{a\l})=\left\{
\begin{array}{ll}
 v^{\mp a^2}  &
{\rm for}\;\; a=2k+1,\;\; k\in\Z \\
 0& {\rm otherwise} \end{array} \right.$$
$$L^1_{\pm 2}(q^{a\l})=\left\{
\begin{array}{ll}
 v^{\mp a^2}  &
{\rm for}\;\; a=2k,\;\; k\in\Z \\
 0& {\rm otherwise} \end{array} \right.$$
This proves b) for $r=4p$. The other case  is similar.

c) From b) we have
$$(L^1_{\pm2}-L^0_{\pm2})(q^{a\l})=(-1)^a (c_1-c_0)L_{\pm2}=(-1)^a(c_1-c_0)v^{\mp a^2}=$$
$$(c_1-c_0)(-v)^{\mp a^2}=(c_1-c_0)L_{\pm2}|_{v\to -v}\, .$$
\end{proof}

\subsection{Cohomological classes}

If $r=2\mod 4$, then the  formula (\ref{def_spin}) defines
cohomological refinements of the  quantum invariant. Here
$\sigma_0$ is the nontrivial cohomological class and $\sigma_1$ is the
other one. We assume throughout
 this subsection that $r=2p$, with odd $p$. Then
 $v^{p^2}=\zeta
 I$, where $I$ is the complex unit and $\zeta=\pm 1$.

\begin{Thm}  Suppose $v$ is a $4p$--root of unity
with $p$ odd and
$v^{p^2}=\zeta I$.

\begin{equation} \tau'_{M_\pm,\sigma_\e}(v) =
\frac{1}{2}\left[ \tau'_M(v) + (-1)^{\e} \zeta  I
\;\frac{1+v}{1-v}\;\;\tilde\tau_M(-v) \right]\label{equ.01}
\end{equation}
\end{Thm}

\begin{proof}
The proof is analogous to the proof of Theorem \ref{four}, replacing Lemma
\ref{as} with Lemma \ref{as2}.
\end{proof}

 It's easy to
see that the right hand side of (\ref{equ.01}) belongs to
$\frac{1}{1-v}\ZZZ$,
and define an invariant of $3$--manifold $M\in \cM_2$ with a fixed
homological structure. This proves the part Theorem
\ref{main_spin} concerning cohomological structures.

 \begin{lem}\label{as2}
$$a) \;\;\;\;
\g^\e_{2,r}=c_\e \g_{2,r}, \;\;\;\;
\g^\e_{-2,r}=c_{\e+1} \g_{-2,r},\;\;\; c_0+c_1=1; $$

$$b)\;\;\;\;\; L_{\pm 2}^0(q^{a\lambda})= \frac{1 \mp \zeta (-1)^a I}{2}\,
 L_{\pm 2} (q^{a\lambda}),$$
$$ L_{\pm 2}^1(q^\lambda a)= \frac{1 \pm \zeta (-1)^a I}{2}\,
 L_{\pm 2}(q^{a\lambda });$$

$$c) \;\;\;
L_{\pm 2}^1(q^{a\lambda }) - L_{\pm 2}^0(q^{a\lambda }) = \pm \zeta
 IL_{\pm 2}(q^{a\lambda})|_{v\to -v}.$$

\end{lem}

 \begin{proof}
$a)$ By shifting $\l\to \l+r/2$, we see that $\gamma^1_{\pm 2,r}=
\pm \zeta I\gamma^0_{\pm2,r}$. This shows that
$c_0=(1-\zeta I)/2$ and $c_1=(1+\zeta I)/2$.

$b)$ From the definition of the odd and even Laplace
transforms we have
$$L^\e_{\pm 2}(q^{a\l})=
\frac{1}{\g_{\pm2,r}} \sum_{\l=\e \pmod 2}q^{\pm\frac{(\l^2-1)}{2}}
q^{a\l}\, .$$
If $a$ is even,
$$L^\e_{\pm 2}(q^{a\l})=\frac{\g^\e_{\pm2,r}}{\g_{\pm2,r}}\;\;
 v^{\mp a^2}. $$
For odd  $a$,
$$L^\e_{\pm 2}(q^{a\l})=
\frac{\g^{\e+1}_{\pm2,r}}{\g_{\pm2,r}}\;\; v^{\mp a^2}. $$
This implies the result.

$c)$ Follows from $b)$ analogously to $c)$ Lemma \ref{as}.

\end{proof}

\section{Quantum invariants for $\Q$HS with $H_1=\Z/b\Z$}

For $M\in \cM_b$, we show that the Laplace transform method
 applies in the case, when
 $gcd(b,r)$ is a power of two.
This leads to a
 construction of  universal invariants
dominating quantum invariants and their refinements
at roots of unity of order $r$
with  $gcd(b,r)=2^n$. As a application, we
 derive    new integrality properties of quantum invariants.


\subsection{Laplace transforms}
We define the odd and even Laplace transforms as follows.
\be\label{lap} \ev_r(L^\e_b(F_n(q^\l,q))):=
\frac{1}{\gamma_{b,r}}\sum^{r-1}_{\l =\e\pmod 2}
q^{\frac{b(\l^2-1)}{4}}\, F_n(q^\l,q)
\,  .\ee

\begin{pro}\label{ref-lap}
For $b=2^t c$ and $r=2^s d$ with $(c,d)=1$,  and $t\neq s+1$,
the odd and even Laplace
transforms are well--defined.
\end{pro}

\begin{proof}
The proof is  by case by case checking.

First assume $s\geq t+2$. Then, by shifting $\l\to \l+2^{s-t}d$,
we see that $\gamma^1_{b,r}=0$.
Furthermore,
$$\sum^{r-1}_{\l=1\;{\rm odd}} q^{b(\l^2-1)/4}\, q^{a\l}=\left\{
\begin{array}{ll}
 q^{-a^2/b} \gamma_{b,r} &
{\rm for}\;\; a=2^{t-1}(2k+1),\;\; k\in\Z \\
 0& {\rm otherwise} \end{array} \right.$$

Indeed,
$$\sum q^{b(\l^2-1)/4}\, q^{a\l}=
q^{-a^2/b} q^{-b/4} \sum q^{(b\l+2a)^2/4b}\, .$$
For $a=2^la'$, $0\leq l < t-1$, $a'$  odd, and $l=t-1$, $a'$ even,
it is easy to see that
the summands for $\l$ and $\l+r/2^{l+1}$ cancel with each other.
For $a=2^{t-1}(2k+1)$,
the sum  is equal to $q^{-a^2/b}\gamma_{b,r}$.

Analogously,     the following formulas
define the even Laplace transform.
$$\sum^{2r-1}_{\l=0\;{\rm even}} q^{b(\l^2-1)/4}\, q^{a\l}=\left\{
\begin{array}{ll}
 q^{-a^2/b} \gamma_{b,r} &
{\rm for}\;\; a=2^{t}k,\;\; k\in\Z \\
 0& {\rm otherwise} \end{array} \right.$$
\v8

For $s=t+1$, we see that $\gamma^0_{b,r}=0$ by shifting $\l\to \l+2d$.
Moreover,
$$L^0_b(q^{a\l})=\left\{
\begin{array}{ll}
 q^{-a^2/b}  &
{\rm for}\;\; a=2^{t-1}(2k+1),\;\; k\in\Z \\
 0& {\rm otherwise} \end{array} \right.$$
$$L^1_b(q^{a\l})=\left\{
\begin{array}{ll}
 q^{-a^2/b}  &
{\rm for}\;\; a=2^{t}k,\;\; k\in\Z \\
 0& {\rm otherwise} \end{array} \right.$$
\v8

If $s=t$, $\gamma^1_{b,r}=\pm I\gamma^0_{b,r}$
(by shifting $\l\to \l+d$) and
$$L^0_b(q^{a\l})=\left\{
\begin{array}{ll}
 (1\pm I)/2\;\, q^{-a^2/b}  &
{\rm for}\;\; a=2^{t-1}(2k+1),\;\; k\in\Z \\
 (1\mp I)/2\;\, q^{-a^2/b}  &
{\rm for}\;\; a=2^{t}k,\;\; k\in\Z \\
 0& {\rm otherwise} \end{array} \right.$$
$$L^1_b(q^{a\l})=\left\{
\begin{array}{ll}
(1\mp I/)2\;\,  q^{-a^2/b}  &
{\rm for}\;\; a=2^{t-1}(2k+1),\;\; k\in\Z \\
(1\pm I)/2\;\, q^{-a^2/b}  &
{\rm for}\;\; a=2^{t}k,\;\; k\in\Z \\
 0& {\rm otherwise} \end{array} \right.$$
\v8
Finally, for $s\leq t-2$, $\gamma^0_{b,r}=\gamma^1_{b,r}$ and
$$L^0_b(q^{a\l})=L^1_b(q^{a\l})
=\left\{
\begin{array}{ll}
1/2\;  q^{-a^2/b}  &
{\rm for}\;\; a=2^s k,\;\; k\in\Z \\
 0& {\rm otherwise} \end{array} \right.$$
\end{proof}

\noindent
{\bf Note.}  If  $t=s+1$, we have  $\gamma^0_{b,r}=-\gamma^1_{b,r}$ and
 $\gamma_{b,r}=0$. Hence, the Laplace transform
cannot be defined by (\ref{lap}) in this case.
But at least if $d=1$,
the method applies, if we normalize the Laplace transform by
$\g^0_{b,r}$ instead of $\g_{b,r}$.
\v8


\subsection{Refined universal invariants}
Let $M\in \cM_b$. Let  $gcd(b,r)=2^n$, $n\in \N$, i.e.
we  assume $b=2^t c$ and $r=2^s d$ with odd $c,d$, and
$gcd(c,d)$=1. If $t\neq s+1$, and $s>0$, then
$$\tau_{L(b,1)}=\frac{L_b(F_0)\g_{b,r}}{2(1-q^{-sn(b)}) \g_{sn(b),r}}
$$ is nonzero and we
can renormalize
$$\tau'_M=\frac{\tau_M}{\tau_{L(b,1)}}\, , \;\;\;\;
\tau'_{M,\sigma_\e}=\frac{\tau_{M,\sigma_\e}}{\tau_{L(b,1)}}\, .$$
If $s=0$,  $\tau^{SO(3)}_{L(b,1)}$ is always nonzero.
In this case, we put
$$\tau'_M=\tau^{SO(3)}_M/{\tau^{SO(3)}_{L(b,1)}}\, .$$

Without loss of generality, we assume that $M=S^3(L)$,
where $L$ is an algebraically split link of $(l+1)$ components,
the  framing of the $0$--th component is $b$, the next $p$ components
have framing $1$, and the remaining ones have framing $-1$.

\begin{thm}
Suppose $M\in \cM_b$, and $b=2^t c$, $r=2^s d$ are as above ($t\neq s+1$).
If $s\neq 0$, the refined quantum invariant of $(M,
\sigma_{\e})$
 is given by the following formula
$$\tau'_{M,\sigma_\e}(q)=
\sum_{k=0}^\infty   \left(\sum_{\max
k_i=k} C_{L,\fk}(v)\,
 \prod\limits^{p}_{i=1}
(-1)^{k_i} q^{-\frac{k_i(k_i+3)}{2}}\frac{L^\e_b(F_{k_0})}{(q^{k_0+1})_{k_0+1}}
\right)\frac{(q^{k+1})_{k+1}}{L_b(F_0)}\, $$
where  $L^\e_b$ are defined in the
proof of Proposition \ref{ref-lap}.
Here
  $\sigma_\e\in H^1(M,\Z/2\Z)$ if $s=1$,  otherwise  $\sigma_\e$
is a spin structure.
If $s=0$,
$$\tau'_{M}(q)=
\sum_{k=0}^\infty   \left(\sum_{\max
k_i=k} C_{L,\fk}(v)\,
 \prod\limits^{p}_{i=1}
(-1)^{k_i} q^{-\frac{k_i(k_i+3)}{2}}\frac{L^1_b(F_{k_0})}{(q^{k_0+1})_{k_0+1}}
\right)\frac{(q^{k+1})_{k+1}}{L_b(F_0)}\, .$$

\end{thm}

\v8
\noindent
{\bf Example.} Suppose $L$ is the Whitehead link with framings $-1$
and $-4$.
Let $M=S^3(L)$.
 Then
$$\tau'_{M,\sigma_\e}(q)=\frac{1}{2}\sum^\infty_{k=0} (-1)^{k} q^{-k(k+1)/2}
L^\e_{-4}(F_k)  \, .$$

\subsection{Proof of Theorem \ref{six}}
Let us first assume that $s>0$.
Then $\tau'_M=\tau'_{M,\sigma_0}+\tau'_{M,\sigma_1}$.
If $t>1$,  $L_b$ sends  $q^{a\l}$ to zero if $a\neq 2^{t-1}k$ with $k\in \Z$,
e.g. if $a$ is odd. We deduce that $L_b(F_0)=2$.
But $L^\e_b$ is  divisible by 2. The result follows.

If $t=1$, $L_b$ sends $q^{a\l}$ to $q^{-a^2/b}$ for all $a$.
Then $L_b(F_0)=2(1-x)$ with $x^{-b}=q$.
We claim that $L_b(F_k)$ is divisible by $L_b(F_0)$ for all $k\in N$.
Indeed, $L_b(F_k)$ can be considered as a polynomial in $x$.
By the $q$--binomial formula, we have
$$L_b(F_k)=2(-1)^k x^{-bk(k+1)/2} \sum^{2k+1}_{j=0}
(-1)^j \left[\begin{array}{c}
2k+1\\ j\end{array}\right]_q x^{\frac{-bj(j-1)}{2}}x^{bkj}x^{(j-k)^2}\, .$$
Then
$$\lim_{x\to 1} L_b(F_k)=2(-1)^k\sum^{2k+1}_{j=0}
(-1)^j \left[\begin{array}{c}
2k+1\\ j\end{array}\right] =0\, .$$

If $s=0$ or $r$ is odd, then $gcd(b,r)=1$ and $L^1_b$ also
sends $q^{a\l}$ to $q^{-a^2/b}$ for all $a$.
The same argument shows that $L^1_b(F_n)$ is  divisible by $2(1-x)$.

$\hfill\Box$

\section*{Appendix }

\setcounter{section}{0}


Here we deduce  Proposition 2.3
 from the Habiro's results in \cite{Ha}.
Let $L$ be an algebraically split  link of $l$ components with all
framings zero.
Let $\fn=\{n_1, n_2,...,n_{l}\}$ be a coloring of $L$ by
$\fn$--dimensional  irreducible
representations of ${\mathfrak sl}_2$.

\v8
\noindent
{\bf Proposition A.} {\it
There exist $C_{L,\fk}(v)\in \Z[v^{\pm 1}]$, such that
$$J'_L(\fn)=\sum_{k=0}^\infty   \left(\sum_{\max k_i=k} C_{L,\fk}(v) \,
(1-q)^l  \prod^{l}_{i=1}
\frac{(q^{1+n_i})_{k_i}(q^{1-n_i})_{k_i}}{(q^{k_i+1})_{k_i+1}}\right)
\frac{(q^{k+1})_{k+1}}{ (1-q)}
$$
}
\v8

\begin{proof}
For a $0$--framed link, we have
$$J_L(\fn)= (-1)^{l-\sum n_i} \la L(e_{\fn-1})\ra\, ,$$
where $\la L(e_{\fn-1})\ra$ is the Kauffman bracket
of $L$, where each component is cabled by $e_{n_i-1}$ (see \cite{B}).
Recall that $\{e_i\}_{i\geq 0}$
provides a basis for the skein algebra of a solid torus.
An other basis is given by elements  $\{R_i\}_{i\geq 0}$
$$R_k=\prod^{k-1}_{i=0}(z-\l_{2i})\, ,\;\;\;\;\l_i=-v^{i+1}-v^{-i-1}\, ,$$
where $z$ is the $0$--framed closed line $S^1\times {pt}$ in the interior
of $S^1\times D^2$, and $z^i$ means $i$ parallel copies of $z$.

The basis change  is given by the following
formula (compare with \cite{Ma})
\be\label{change}
e_{n-1}=\sum^{n-1}_{k=0} (-1)^{n-1-k}\left[\begin{array}{l}
n+k\\n-1-k\end{array}\right] R_k\, ,\;\;\;{\rm where}\ee
$$\left[\begin{array}{l}
a\\b\end{array}\right]=\frac{[a]!}{[b]![a-b]!} \;\;\;\;[a]!=\prod^a_{i=1}
[i]\, .$$
Using (\ref{change}) we get
$$J_L(\fn)=\sum_{\fk} (-1)^{l-\sum_i ( n_i-k_i)}\prod^{l}_{i=1}
\left[\begin{array}{l}
n_i+k_i\\n_i-1-k_i\end{array}\right]
J_L(R_{k_0}, R_{k_1},...,R_{k_{l-1}})\, .$$
The crucial step in the proof is Theorem 3.3 in \cite{Ha}.
The first part of Theorem 3.3  provides
the existence of $c_{L,\fk}\in \Z[v^{\pm 1}]$
such that
$$J'_L(\fn)=\sum_{\fk}c_{L,\fk} (1-q)^l \prod^{l}_{i=1}
\; \frac{S(n_i,k_i)}{(q^{k_i+1})_{k_i+1}}\, ,$$
where
$$S(n,k)=\frac{\{n-k\}\{n-k+1\}...\{n+k\}}{\{n\}}\, ,\;\;\;\; \{i\}=
v^i-v^{-i}\, .$$
The second part of  Theorem 3.3 implies the result.
\end{proof}

\bibliographystyle{amsplain}

\end{document}